\documentclass[12pt]{article}
\usepackage{amsmath, amssymb,color}
\usepackage{hyperref}

\newcommand{\nc}{\newcommand}
\nc{\ga}{\gamma}
\nc{\di}{\displaystyle}
\nc{\ek}{\protect\\[1ex]}
\nc{\C}{{\mathbb C}}
\nc{\N}{{\mathbb N}}
\nc{\R}{{\mathbb R}}
\nc{\Z}{{\mathbb Z}}
\nc{\La}{\Lambda} \nc{\la}{\lambda}
\nc{\da}{\delta}
\nc{\Da}{\Delta}
\nc{\na}{\nabla}
\nc{\vp}{\varphi}
\nc{\ka}{\kappa}
\nc{\si}{\sigma}
\nc{\Si}{\Sigma}
\nc{\al}{\alpha}
\nc{\be}{\beta}
\nc{\om}{\omega}
\nc{\Om}{\Omega}
\nc{\pa}{\partial}
\nc{\ti}{\times}
\nc{\ve}{\varepsilon}
\nc{\ra}{\rightarrow}
\nc{\Ra}{\Rightarrow}
\nc{\ran}{\rangle}
\nc{\lan}{\langle}
\nc{\eq}[1]{\mbox{\rm{(\ref{E#1})}}}
\nc{\EL}{\cal L}
\nc{\qed}{\mbox{}\nolinebreak\hfill \rule{2mm}{2mm}}
\nc{\ha}{\frac{1}{2}}
\nc{\wrar}{\rightharpoonup}
\nc{\hra}{\hookrightarrow}
\nc{\supp}{\text{supp}\,}
\nc{\curl}{\text{curl}\,}
\nc{\dense}{\hra^{\hspace{-3mm}d\,\,}}

\newtheorem{lem}{Lemma}[section]
\newtheorem{theo}[lem]{Theorem}
\newtheorem{coro}[lem]{Corollary}
\newtheorem{propo}[lem]{Proposition}
\newtheorem{rem}[lem]{Remark}

\renewcommand{\div}{{\rm{div}}\,}

\numberwithin{equation}{section} 

\begin{document}
\bibliographystyle{alpha}
\title{Maximal $L^1$-regularity of generators for bounded analytic semigroups in Banach spaces}
\author{Myong-Hwan Ri\thanks{Institute of Mathematics, State Academy of
Sciences, DPR Korea; email:\vspace{-3pt}
\texttt{math.inst@star-co.net.kp}} \;and
Reinhard Farwig\thanks{Department of Mathematics, Darmstadt University of
Technology, Germany; email:\vspace{-3pt}  \texttt{  farwig@mathematik.tu-darmstadt.de}}
}

\date{ }

\maketitle

\begin{abstract}
\noindent In this paper, we prove that the generator of any bounded
analytic semigroup in $(\theta,1)$-type real interpolation of its
domain and underlying Banach space has maximal $L^1$-regularity,
using a duality argument combined with the result of maximal
continuous regularity.
 As an application, we consider maximal
$L^1$-regularity of the Dirichlet-Laplacian and the Stokes operator
in inhomogeneous $B^s_{q,1}$-type Besov spaces on domains of
$\mathbb R^n$, $n\geq 2$.
\end{abstract}
{\small{\bf 2000 Mathematical Subject Classification:} 35K90; 46B70; 47D06



\section{Introduction}

A central question for parabolic evolution equations in Banach
spaces is whether a related linear operator will have {\it maximal
regularity}. Maximal regularity of a closed linear operator can
efficiently be applied to study related quasilinear and semilinear
evolution equations; for relevant background material we refer {\em
e.g.} to \cite{Am95}, \cite{DHP03} and \cite{Pr02}. In this paper we
study maximal $L^1$-regularity of sectorial operators in Banach
spaces. An advantage of maximal $L^1$-regularity of sectorial
operators is that for the Cauchy problem \eq{2.8} with nonzero
initial value $u_0$ the existence of strong solutions is guaranteed
without requiring higher regularity on $u_0$ than required for the
mild solution; this is not the case with maximal $L^p$-regularity
for $p\in (1,\infty)$. This property is efficient, for example, when
considering density-dependent Navier-Stokes equations, {\em cf.}
\cite{DaMu09}, \cite{DaMu15}.

Given a Cauchy problem
\begin{equation}
\label{E2.8} u'(t) + Au(t)= f(t)\quad\text{in } X,\quad u(0) = 0,
\end{equation}
 with a closed linear operator $A$ in a complex Banach space $X$,
the operator $A$ is said to have {\it maximal $L^p$-regularity} in
$X$ for $p\in [1,\infty]$, if for every $f\in L^p(\R_+;X)$, where
$\R_+=[0,\infty)$, the problem \eq{2.8} has a unique solution
satisfying
\begin{equation}
\label{E2.9} \|u_t\|_{L^p(\R_+;\;X)}+\|Au\|_{L^p(\R_+;\;X)}\leq
C\|f\|_{L^p(\R+;X)}\quad (\exists C>0).
\end{equation}

Let $J=[0,T]$ for $T<\infty$ or $J=[0,\infty)$. The operator $A$ in
$X$ is said to have {\it maximal continuous regularity} on $J$ if
for every $f$ belonging to $BC(J,X)$ the problem \eq{2.8} has a
unique solution satisfying
\begin{equation}
\label{E2.10} \|u_t\|_{BC(J,X)}+\|Au\|_{BC(J,X)}\leq
C\|f\|_{BC(J,X)} \quad(\exists C>0),
\end{equation}
where $BC(J,X)$ is  the space of all $X$-valued bounded and
continuous functions on $J$.

 Let us recall some known results for maximal regularity. If a
sectorial operator in a Banach space has maximal $L^p$-regularity
for some $p\in (1,\infty)$, then so does it for all $p\in
(1,\infty)$ (see \cite{Can86}, \cite{Do00}, \cite{Sob64});
furthermore, \eq{2.9} holds with temporally weighted space
$L^p_\ga(\R_+;X)$ with norm $\|t^{1-\ga}u\|_{L^p(\R_+;X)}$ for all
$p\in (1,\infty)$, $\ga\in (1/p,1]$, {\em i.e.} {\it maximal
$L^p_\ga$-regularity}, see \cite{PrSi04}. If a densely defined
closed linear operator $A$ has maximal $L^p$-regularity for $p\in
(1,\infty)$ or maximal continuous regularity, then $-A$ must
generate an analytic semigroup in the underlying Banach space
(\cite{Do93}, \cite{Do00}, \cite{LeSi11}).

By a classical result (\cite{Sim64}) generators of bounded analytic
semigroups in Hilbert spaces have maximal $L^p$-regularity for $p\in
(1,\infty)$. It was shown in \cite[Theorem 3.2]{DoVe87} that, if a
sectorial operator $A$ in a {\em UMD} space $X$ admits bounded
imaginary powers ({\em BIP}) in ${\cal L}(X)$ with power angle less
than $\frac{\pi}{2}$, {\em i.e.}, $\|A^{i\tau}\|_{{\cal L}(X)}\leq
ce^{\theta |\tau|}$, $\tau\in\R$, with some $\theta\in
(0,\frac{\pi}{2})$, then $A$ has maximal $L^p$-regularity in $X$ for
$p\in (1,\infty)$, see also \cite{CDMY96}, \cite[Theorem
4.4]{DHP03}, \cite[Theorem III.4.10.7]{Am95}. Note that {\em UMD}
spaces are reflexive.

The operator $A$ with $-A$ being the generator of a bounded analytic
semigroup in a Banach space $E$ has maximal continuous and
$L^\infty_\ga (\ga\in (0,1])$- regularity in the continuous
interpolation space $X=(E,E_1)^0_{\theta,\infty}$, $\theta\in
(0,1)$, where $E_1=\mathcal D(A)$ is endowed with the graph norm of
$A$, see \cite[Th\'eor\`eme 3.1]{DaGr79}, \cite[Theorem 2.14]{An90},
\cite[Theorem III.3.4.1]{Am95} and \cite[Theorem 3.5]{RiZhZh15}. We
note that if $A$ has maximal continuous regularity in a Banach space
$X$, then $X$ can be neither reflexive nor weakly sequentially
complete (\cite[Remark 2.4 (d)]{ClSi01}, \cite[Remarks III.3.1.3
(b)]{Am95}). Moreover, if $A$ is unbounded and has maximal
$L^\infty$-regularity, then $X$ must contain a complemented subspace
isomorphic to $c_0$ (\cite{Gu95}).


\par\medskip
Concerning maximal $L^1$-regularity, in \cite{GiSa12}, the property
was proved for the negative Laplacian in $FM(\R^n)$, which is the
Fourier image of the space of Radon measures on $\R^n$ with finite
total variation. It is shown that the Stokes operator has maximal
$L^1$-regularity in homogeneous Besov spaces
$\dot{B}^s_{p,1}(\Omega)$, $p\in (1,\infty)$, $s\in (1/p-1,1/p)$,
where $\Omega$ is the whole or half space, a smooth bounded or
exterior domain, see \cite{DaMu09},\cite{DaMu15}. It is worth
mentioning that in these papers the proof is essentially based on
properties of the heat kernel on $\R^n$; moreover, in \cite{DaMu09},
\cite{DaMu15} the Littlewood-Paley characterization of Besov spaces
is crucial.

On the other hand, for generator $-A$ of an analytic semigroup, $A$
does not have maximal $L^1$-regularity property if it is unbounded
and the underlying Banach space does not contain a complemented
subspace isomorphic to $l_1$, see \cite[Theorem 5]{Gu95}. Thus, it
follows by \cite[Theorem 2.e.7]{LiTz77} that maximal
$L^1$-regularity does not hold in reflexive Banach spaces. It is
known by \cite[Theorem 3.6]{KaPo08} that when $-A$ is the generator
of a bounded analytic semigroup in a Banach space $X$, maximal
$L^1$-regularity of $A$ is equivalent to
\begin{equation}
\label{E1.4} \int_0^\infty \|Ae^{-tA}u\|_Xdt \leq C\|u\|_X,\forall
u\in X.
 \end{equation}
One can check that \eq{1.4} holds when $A$ is the generator of a
bounded analytic semigroup in a Banach space $E$  with $0$ in its
resolvent set $\rho(A)$ and $X=(E,E_1)_{\theta,1}$, $\theta\in
(0,1)$, where $E_1={\cal D}(A)$ is endowed with graph norm of $A$.
However, it is not clear whether \eq{1.4} can still be easily proved
without the invertibility condition $0\in\rho(A)$. We note that
maximal regularity property without assuming $0\in\rho(A)$ allows
the estimate constant to be independent of time interval, which is
very important in existence and stability theory for corresponding
nonlinear problems.

\par\medskip
In this paper we show that maximal $L^1$-regularity for the
generator of any bounded analytic semigroup in $(\theta,1)$-type
real interpolation of Banach spaces can follow directly from the
{\it known} maximal continuous regularity result by a duality
argument without assuming $0\in\rho(A)$. We remark here that
$(L^1(J,E))'=L^\infty(J,E')$ and $(BC(J,E))'={\cal M}(J,E')$, the
space of all $E'$-valued BV measures, and hence the duality argument
is not obvious at all.

Below, for Banach spaces $E$ and $E_1$, the notation $E_1\dense E$
means that $E_1$ is continuously and densely embedded in $E$. The
dual space of $E$ is denoted by $E'$. Let ${\cal L}(E_1,E)$ stand
for the set of all linear bounded operators from $E_1$ to $E$, and
let ${\cal L}(E)={\cal L}(E,E)$. Furthermore, $\text{Isom}(E_1,E)$
denotes the set of all norm isomorphisms from $E_1$ to $E$. The
notation $A\in {\cal H}(E_1, E)$ means that $A\in {\cal L}(E_1, E)$
and, considered as an unbounded linear operator in $E$ with
$\mathcal D(A)=E_1$, $-A$ generates an analytic $C_0$-semigroup
$\{e^{-tA}\}_{t\geq 0}$ in $E$.

The main result of the paper is as follows.

\begin{theo} [Maximal $L^1$-regularity] \label{T1.1}
Let $E$ and $E_1$ be complex reflexive Banach spaces  with
$E_1\dense E$ and let $A\in {\cal H}(E_1,E)$ generate a bounded
analytic semigroup in $E$. Then, for $\theta\in (0,1)$ the
realization of $A$ in $E_{\theta,1}:=(E,E_1)_{\theta,1}$ has maximal
$L^1$-regularity in $E_{\theta,1}$. More precisely, for $f\in
L^1(\R_+;E_{\theta,1})$, the Cauchy problem \eq{2.8} with
$X=E_{\theta,1}$ has a unique solution $u$ such that $u_t,\, Au \in
L^1(\R_+; E_{\theta,1})$ and
 \begin{equation}
  \label{E2.20}
\begin{array}{rl}
 \|u_t\|_{L^1(\R_+; E_{\theta,1})} + \|Au\|_{L^1(\R_+; E_{\theta,1})}  &
 \leq \;c\|f\|_{L^1(\R_+; E_{\theta,1})},
\end{array}
\end{equation}
 where $c=c(\theta)$ is independent of $f$.

\end{theo}

The main idea for the proof of Theorem \ref{T1.1} is to use a duality
argument based on maximal continuous regularity for a backward-in-time
problem in continuous interpolations spaces which are obtained by a
method of extrapolation of $E_1$ and $E$ introduced {\em e.g.} in
\cite{Am95}, \cite{Am00} (see also Subsection 2.1).

The following corollary is a direct consequence of Theorem \ref{T1.1}.

\begin{coro}
\label{C1.2}  Under the same assumptions on the spaces $E, E_1$
and the sectorial operator $A$ as in Theorem \ref{T1.1}, let $\{(E_\al, A_\al):\al\in\R\}$ be
the interpolation and extrapolation scale generated by the real
interpolation functor $\{(\cdot,\cdot)_{\theta,1}: \theta\in
(0,1)\}$ and $(E,A)$, see Subsection 2.1. Moreover, if $\al\in \Z$, let
$$ E^\bullet_\al:=(E_{\al-1/2}, E_{\al+1/2})_{1/2,1}. $$

(i) If $\al\in (k,k+1)$ with $k\in\Z$, then $A_\al$ has maximal
$L^1$-regularity in $E_{\al-1}$.

(ii) If $\al\in \Z$, then $A^\bullet_\al$, the
$E^\bullet_\al$-realization of $A_{\al-1/2}$,
has maximal $L^1$-regularity in $E^\bullet_\al$.
\end{coro}


This paper is organized as follows. In Section 2, preliminaries on
interpolation and extrapolation scales and vector measures with
bounded variation are given. Section 3 is devoted to the proof of
main results. As an application of the abstract theory, in Section
4, maximal $L^1$-regularity of the Dirichlet-Laplacian and the
Stokes operator in $B_{q,1}^s$-type Besov spaces on domains is
considered.

\section{Preliminaries}
\vspace*{-3mm}
\subsection{Interpolation and extrapolation scales}
Given an interpolation couple $(E_0,E_1)$ of Banach spaces and
$\theta\in (0,1)$, we denote by $[\cdot,\cdot]_\theta$,
$(\cdot,\cdot)_{\theta,r}$, $1\leq r\leq\infty$, the complex and
real interpolation functors, respectively, {\em cf.} {\cite{BL77},
\cite{Tri83}. The {\it continuous interpolation space}
$(E_0,E_1)^0_{\theta,\infty}$, $\theta\in (0,1)$, of $E_0$ and $E_1$
is defined as the closure of $E_0\cap E_1$ in
$(E_0,E_1)_{\theta,\infty}$. }

To be more specific, in this subsection we assume that  $E_1\dense E$ and $E$ is reflexive. Furthermore, let an interpolation functor
$(\cdot,\cdot)_\theta$, $\theta\in (0,1)$, be given by
\begin{equation}
\label{E2.5} (\cdot,\cdot)_\theta\in
\big\{[\cdot,\cdot]_\theta,(\cdot,\cdot)_{\theta,p},
(\cdot,\cdot)^0_{\theta,\infty}: 1\leq p<\infty\big\}.
\end{equation}

Let $A$ be a linear, closed and sectorial operator
 (with angle $\vartheta\in (0,\pi)$) in $E_0$
 with $\mathcal D(A)=E_1$ and $\overline{R(A)}=E$, {\em i.e.,}
 $S_\vartheta := \{\lambda\in \C\setminus\{0\}: |\arg \lambda|<\vartheta\}\subset \rho(-A)$ and
\begin{equation}
\label{E2.2n} \|\lambda(\lambda+A)^{-1}\|_{{\cal L}(E_0)}\leq
K\quad\textrm{ for all }\; \lambda \in S_\vartheta
 \end{equation}
  with a constant
$K>0$ depending on $\vartheta$. Suppose that $A\in {\cal
H}(E_1, E)$ and
\begin{equation}
\label{E2.3n} \|e^{-tA}\|_{{\cal L}(E)}\leq M \quad (\exists M>0,
\forall t\geq 0).
\end{equation}

For $k\in \N$ let  $E_k:=\mathcal D(A^k)$ be the domain of $A^k$
in $E_0:=E$ endowed with its graph norm. Let $A^\natural:=A'$ denote
the dual (adjoint) operator of $A$ in the dual space $E^\natural=E'$
which is a closed, densely defined unbounded operator with domain
$\mathcal D(A^\natural) \subset E^\natural$. Then, for {\em
negative} integers $k$, the space $E_{k}$ can be introduced as the
dual space of $\mathcal D((A^\natural)^{-k})$ which is the domain of
$(A^\natural)^{-k}$ in $E^\natural$ endowed with its graph norm.

Finally, for $\al\in (k,k+1),\; k\in\Z$, the space  $E_\al$ is defined by
\begin{equation} \label{E2.4n}
E_\al:= (E_k, E_{k+1})_{\al-k},
\end{equation}
and the operator $A_\al$ in $E_\al$ by the realization of $A$ in
$E_\al$ if $\al\geq 0$ and as the closure of  $A$ in $E_\al$ if
$\al<0$.
Thus, the {\em interpolation and extrapolation scale} $\{(E_\al,
A_\al):\al\in\R\}$ is generated by $(\cdot,\cdot)_\theta$ and
$(E,A)$. It is well-known that for $-\infty<\beta<\al<\infty$
\begin{equation}
 \label{E2.6}
E_\al\dense E_\beta,\; A_\al\in {\cal H}(E_{\al+1},
E_\al), \quad
e^{-tA_\al}=e^{-tA_\be}|_{E_\al},\quad \|e^{-tA_\al}\|_{{\cal
L}(E_\al)}\leq M
\end{equation}
with the same constant $M>0$ as in \eq{2.3n}, see \cite{Am88},
\cite[Chapter V.1, V.2]{Am95}, \cite[Section 1]{Am00}.

From now on, suppose that
$(\cdot,\cdot)_\theta$ for every $\theta\in (0,1)$ does not coincide
with $(\cdot,\cdot)^0_{\theta,\infty}$ and let $\{(E^\natural_\al,
A^\natural_\al):\al\in\R\}$ be the interpolation and extrapolation
scale generated by $\{E^\natural, A^\natural\}$ and the
interpolation functor
\begin{equation} \label{E2.4nn}
(\cdot,\cdot)^\natural_\theta=
    \left\{
      \begin{array}{ll}
      [\cdot, \cdot]_\theta & \quad\text{for } (\cdot,\cdot)_\theta=[\cdot,\cdot]_\theta,\ek
     (\cdot, \cdot)_{\theta,r'} & \quad \text{for } (\cdot,\cdot)_\theta=(\cdot,\cdot)_{\theta,r},  1<r<\infty,\ek
    (\cdot, \cdot)^0_{\theta,\infty} & \quad\text{for }
    (\cdot,\cdot)_\theta=(\cdot,\cdot)_{\theta,1}.
     \end{array}
   \right.
\end{equation}
Note that compared to the literature, see {\em e.g.}
\cite[p. 282]{Am95}, \cite[(1.9)]{Am00}, the role of
$(\cdot,\cdot)_\theta$ and $(\cdot,\cdot)^\natural_\theta$ is
interchanged since the third case in \eq{2.4nn} starts with
$(\cdot,\cdot)_{\theta}=(\cdot,\cdot)_{\theta,1}$ and defines the
continuous interpolation functor $(\cdot,\cdot)_{\theta,\infty}^0$.
It is this case which is used in this article. Hence, in view of
reflexivity of $E$ and $(A^\natural)'=A$ due to closedness of the
operator $A$ in $E$, results from the literature can be used by
formally interchanging $E$ and $E^\natural$, {\em i.e.},
interchanging $(E_\al, A_\al)$ with $(E^\natural_\al,
A^\natural_\al)$.

As an example we will frequently use that
  \begin{equation}
  \label{E2.12}
  E_\al=(E^\natural_{-\al})', \quad \al\in\R,
  \end{equation}
{\em cf.} \cite[Theorem V.1.5.12]{Am95}, \cite[(1.10)]{Am00},
whereas $(E_\al)'$ will be different from  $E^\natural_{-\al}$ by
the third case of \eq{2.4nn}. Moreover, let us mention the
smoothness properties of $e^{-tA_\al}$, see
\cite[Theorem 6.1 (vi)]{Am88}, \cite[Theorem
V.2.1.3, Corollary V.2.1.4]{Am95}):
 For $-\infty<\beta<\alpha<\infty$ and $f\in E_{\be}$
 there holds $t^{\al-\be}e^{-tA_{\al}}f\in BC((0,T],E_{\al})$ and
 \begin{equation}
\label{E2.4}
\|e^{-tA_{\al}}f\|_{E_{\al}}
 \leq c(\al,\beta,M)t^{\be-\alpha} \|f\|_{E_{\be}}, \quad t\in(0,T].
\end{equation}

Let $(A^\natural_{-\al-1})'$ be the dual operator of
$A^\natural_{-\al-1}\in {\cal L}(E^\natural_{-\al},
E^\natural_{-\al-1})$, $\al\in\R$, in the sense of bounded linear
operators. In view of \eq{2.12}, $(A^\natural_{-\al-1})'\in {\cal
L}((E_{-\al-1}^\natural)', (E_{-\al}^\natural)') = {\cal
L}(E_{\al+1}, E_\al)$. Furthermore, the above argument of
interchanging $(E_\al, A_\al)$ by $(E^\natural_\al,
A^\natural_\al)$ and \cite[Theorem V.2.3.2]{Am95} imply that
\begin{equation}
\label{E2.13} (A^\natural_{-\al-1})'=A_\al \in {\cal L}(E_{\al+1},
E_\al),\;  \al\in\R.
\end{equation}

\subsection{Vector measures with bounded variation}
Let $(E,\|\cdot\|_E)$ be a Banach space, $J$ be a $\si$-compact
metrizable space, and let ${\cal B} = {\cal B}(J)$ be the Borel
$\si$-algebra of $J$. A $\si$-additive map $\mu: {\cal B}\mapsto E$
is said to be an $E$-valued vector measure if $\mu(\emptyset)=0$
({\em cf.} \cite[Subsection 2.2]{Am03}. For a vector measure $\mu$
the {\it total variation}
 $|\mu|: {\cal B}\mapsto \R_+\cup\{\infty\}$ is defined by
$$|\mu|(G):=\sup_{\pi(G)}\sum_{F\in \pi(G)}\|\mu(F)\|_{E},\quad G\in {\cal B},$$
where the supremum is taken over all partitions $\pi(G)$ of $G$ into
a finite number of pairwise disjoint Borel subsets. Then, $\mu$ is
said to be of {\it bounded variation} on $J$ if
$$\|\mu\|_{BV}:=|\mu|(J)<\infty.$$
We denote by ${\cal M}(J, E):=({\cal M}(J, E),
\|\cdot\|_{BV})$ the space of all $E$-valued vector measures on $J$
with bounded total variation.

Next we replace $E$ in ${\cal M}(J, E)$ by $E'$. Then, through the duality pairing
\begin{equation}
\label{E2.22}
 \lan \mu, u\ran_{{\cal M}(J, E'),BC(J, E)} = \int_J \;
\lan u,\, d\mu\ran_{E,E'}\,,\quad \mu\in {\cal M}(J, E'),\; u\in
BC(J,E),
\end{equation}
it holds by \cite[Theorem 2.2.4]{Am03}
\begin{equation}\label{Edual}
(BC(J, E))'={\cal M}(J, E').
\end{equation}
If  $h\in L^1(J,E')$, then the $E'$-valued measure $\mu_h$ on $J$
defined by
\begin{equation}
\label{E2.23} \mu_h(B):=\int_B h(t)\,dt, \quad  B\in{\cal B},
\end{equation}
has bounded total variation, and the map $h \mapsto \mu_h$
is a linear isometry from  $L^1(J, E')$ into ${\cal
M}(J,E')$. Hence $L^1(J, E')$ can be identified with a closed
subspace of ${\cal M}(J,E')$, {\em cf.} \cite[Remark 2.2.1]{Am03}.

\section{Proof of the main result}

\noindent {\bf Proof of Theorem \ref{T1.1}.}
Let $\{(E_\al, A_\al):\al\in\R\}$ be the
interpolation and extrapolation scale generated by the interpolation
functor $\{(\cdot,\cdot)_{\theta,1}: \theta\in (0,1)\}$ and $(E,A)$,
and let $\{(E^\natural_\al, A^\natural_\al):\al\in\R\}$ be the
interpolation and extrapolation scale generated by the interpolation
functor $\{(\cdot,\cdot)^0_{\theta,\infty}: \theta\in (0,1)\}$ and
$(E^\natural,A^\natural)\equiv (E',A')$, see Subsection 2.1. Then,
by \eq{2.4n} and \eq{2.12}, it holds that
\begin{equation}
\label{E2.17} \begin{array}{rl}
  E_{\al+k}=(E_k, E_{k+1})_{\al,1}, \qquad&
\al\in (0,1),\; k\in\N\cup\{0\},\ek
 E^\natural_{-\al-k}=(E^\natural_{-k},
E^\natural_{-k-1})^0_{\al,\infty},\; & \al\in (0,1), \;k\in\N\cup\{0\},
\end{array}
\end{equation}
and
\begin{equation}
\label{E3.2} E_{\al}=(E^\natural_{-\al})', \quad\al\in\R.
 \end{equation}

Now, let $0<T<\infty$ and consider the backward-in-time Cauchy
problem
\begin{equation}
\label{E2.2} -v_t+A^\natural_{-1-\theta} v=A^\natural_{-1-\theta}
g\quad\text{for } 0\leq t<T,\quad v(T)=0,
\end{equation}
where $g\in C([0,T], E^\natural_{-\theta})$. Note that \eq{2.2} is
reduced by the change of variable $t\mapsto T-t$ to a parabolic
Cauchy problem in $E^\natural_{-1-\theta}$ with initial time $t=0$
and has a unique mild solution $v\in
C([0,T],E^\natural_{-1-\theta})$ expressed by
$$v(T-t)=\int_0^t
           e^{-(t-\tau)A^\natural_{-1-\theta}}A^\natural_{-1-\theta}g(T-\tau)\,d\tau,\;t\in (0,T).$$
 By the
 well-known property of analytic semigroups and the fact that
 $A^\natural_{-1-\theta}$ is an extension of $A^\natural_{-\theta}$, we have
$$\begin{array}{rcl}
v(T-t)&=& A^\natural_{-1-\theta} \int_0^t
           e^{-(t-\tau)A^\natural_{-1-\theta}}g(T-\tau)\,d\tau\ek
       &=& A^\natural_{-\theta} \int_0^t
           e^{-(t-\tau)A^\natural_{-1-\theta}}g(T-\tau)\,d\tau\ek
      &=& A^\natural_{-\theta} w(T-t),\;t\in (0,T),
\end{array}$$
where $w$ is the (unique) mild solution to the backward problem
$$-w_t+A^\natural_{-\theta} w=g\quad\text{in }(0,T),\quad w(T)=0.$$
 Since the
operator $A^\natural_{-\theta}$ is the
$E^\natural_{-\theta}$--realization of  $A^\natural_{-1}\in
{\mathcal H}(E^\natural_{0},E^\natural_{-1})$ and
$E^\natural_{-\theta}=(E^\natural_{0},
E^\natural_{-1})^0_{\theta,\infty}$ by \eq{2.17}$_2$, we get by the
result of maximal continuous regularity (see \cite[Theorem
3.1]{DaGr79}, \cite[Theorem on p.45]{ClSi01}, or \cite[Theorem
III.3.4.1 ($\mu=0$)]{Am95}) that
$$\|A^\natural_{-\theta} w\|_{C([0,T],
E^\natural_{-\theta})} \leq C\|g\|_{C([0,T], E^\natural_{-\theta})}
$$
with constant $C>0$ independent of $T$. Consequently, we have
\begin{equation}
\label{E2.3} \|v\|_{C([0,T], E^\natural_{-\theta})} \leq
C\|g\|_{C([0,T], E^\natural_{-\theta})}.
\end{equation}

Now, let $f\in L^p(0,T; E_1)$, $1<p<1/(1-\theta)$. Then
the mild solution to \eq{2.8} with $X=E_\theta$ ($\equiv
E_{\theta,1}$) is given by
\begin{equation}
\label{E2.1}
 u(t)\equiv u_f(t):=\int_0^t e^{-(t-s)A}f(s)\,ds,
 \end{equation}
and by \eq{2.4} it follows that $\|Ae^{-tA}\|_{{\cal L}(E_1,
E_{\theta})}\leq c\|e^{-tA}\|_{{\cal L}(E_1, E_{1+\theta})} \leq
c/t^{\theta}$. Therefore,
\begin{align*}
\|Au(t)\|_{E_{\theta}}
& \leq  \int_0^t \|Ae^{-(t-s)A}f(s)\|_{E_{\theta}}\,ds\ek
& \leq c\int_{\R^1} \frac{1}{(t-s)^{\theta}}\|\tilde{f}(s)\|_{E_{1}}\,ds, \quad t\in (0,T),
\end{align*}
where  $\tilde{f}$ denotes the extension of $f$ by $0$ from $[0,T]$ to
$\R^1$.
 Hence, by the Hardy-Littlewood-Sobolev inequality we get that
$$\|Au\|_{L^{p*}(0,T;E_{\theta})}\leq C_p\|\tilde{f}\|_{L^p(\R; E_1)}=C_p\|f\|_{L^p(0,T; E_1)}$$
for $p^*=p/(1-(1-\theta) p)$. Moreover, since
$\|e^{-tA}\|_{{\cal L}(E_\theta)}\leq M$, see \eq{2.6}, it follows
from \eq{2.1} that $\|u\|_{L^\infty(0,T;E_\theta)}\leq
M\|f\|_{L^1(0,T;E_\theta)}$ and hence $u\in L^1(0,T;E_\theta)$ in
view of $T<\infty$.

Thus we have
 \begin{equation}
 \label{E2.7}
 Au\in L^{1}(0,T;E_{\theta}),\; u\in L^{1}(0,T;E_{1+\theta})
 \end{equation}
 for $T<\infty$.
Here, recall that $(E^\natural_{-\theta})'=E_\theta$ but
$(E_\theta)'\neq E^\natural_{-\theta}$. Hence
$L^{1}(0,T;E_{\theta})$ is {\em not} the dual space of $L^\infty(0,T;
E^\natural_{-\theta})$ and vice versa.

By the properties of $v, v_t, u, u_t$ in \eqref{E2.3}, \eq{2.1} and
\eq{2.7} it follows by an approximation argument that for almost all
$t\in (0,T)$
$$\frac{d}{dt}\lan u(t), v(t)\ran_{E_\theta, E^\natural_{-\theta}}
 =\lan u_t(t), v(t)\ran_{E_\theta, E^\natural_{-\theta}}
 +\lan u(t), v_t(t)\ran_{E_{1+\theta}, E^\natural_{-1-\theta}}\in L^1(0,T);$$
hence the map $t\mapsto \lan u(t), v(t)\ran_{E_\theta,
E^\natural_{-\theta}}$ is absolutely continuous in $[0,T]$.
Therefore, in view of  $u(0)=v(T)=0$, we have
 \begin{align}
 \label{E3.3}
0 & = \int_0^T \frac{d}{dt}\lan u(t), v(t)\ran_{E_\theta,
E^\natural_{-\theta}}\,dt \nonumber\ek
& = \int_0^T \big(\lan u_t(t), v(t)\ran_{E_\theta, E^\natural_{-\theta}}
    +\lan u(t), v_t(t)\ran_{E_{1+\theta}, E^\natural_{-1-\theta}}\big)\,dt.
\end{align}
Since, by \eq{2.13}, the dual of $A^\natural_{-1-\theta} \in {\cal
L}(E^\natural_{-\theta}, E^\natural_{-1-\theta})$ equals
$A_\theta\in {\cal L}(E_{1+\theta}, E_\theta)$, we get from \eq{2.2}
and \eq{3.3} that
\begin{align*}
 \int_0^T \lan Au, g\ran_{E_\theta,E^\natural_{-\theta}}\,dt
 & = \int_0^T\lan u, A^\natural_{-1-\theta} g\ran_{E_{1+\theta},E^\natural_{-1-\theta}}\,dt\nonumber\ek
 & = \int_0^T\lan u, -v_t+A^\natural_{-1-\theta} v\ran_{E_{1+\theta},E^\natural_{-1-\theta}}\,dt\nonumber\ek
 & = \int_0^T\lan u_t+A_\theta u,
v\ran_{E_{\theta},E^\natural_{-\theta}}\,dt\nonumber\ek
& = \int_0^T \lan f,v\ran_{E_{\theta},E^\natural_{-\theta}}\,dt.
 \end{align*}

In view of \eq{2.22}, \eq{dual}, \eq{2.23} it follows that
 $$
 \begin{array}{l}
 \lan \mu_{Au}, g\ran_{{\cal M}([0,T],E_{\theta}), C([0,T], E^\natural_{-\theta})} = \lan
\mu_f, v\ran_{{\cal M}([0,T],E_{\theta}), C([0,T],
E^\natural_{-\theta})},
 \end{array}
$$
and hence, by \eq{2.3},
\begin{equation*}
\begin{array}{rcl}
 |\lan \mu_{Au}, g\ran_{{\cal M}([0,T],E_{\theta}), C([0,T], E^\natural_{-\theta})}|
 & \leq &
\|\mu_f\|_{{\cal M}([0,T],E_{\theta})}\|v\|_{C([0,T],
E^\natural_{-\theta})}\ek
 & \leq &
C(\theta)\|\mu_f\|_{{\cal M}([0,T],E_{\theta})}\|g\|_{C([0,T],
E^\natural_{-\theta})}.
\end{array}
\end{equation*}
Moreover, since $g$ is an arbitrary element of $C([0,T],E^\natural_{-\theta})$, we have
\begin{equation}
\label{E2.15} \|\mu_{Au}\|_{{\cal M}([0,T],E_{\theta})}\leq
C(\theta)\|\mu_f\|_{{\cal M}([0,T],E_{\theta})}\quad \forall
f\in L^p(0,T; E_1),
\end{equation}
where $C(\theta)$ is independent of $T$. Then by the isometry property of the map $h\mapsto
\mu_h$ we conclude from \eq{2.15} that
\begin{equation}
\label{E2.21}
 \|Au\|_{L^1(0,T;E_\theta)}\leq C(\theta)\|f\|_{L^1(0,T;E_\theta)}\quad
\,\forall f\in L^p(0,T;E_1).
\end{equation}

Now, due to the density of $L^p(0,T; E_1)$ in $L^1(0,T; E_{\theta})$
and linearity of the problem \eq{2.8}, it follows that \eq{2.21}
holds for all $f\in L^1(0,T;E_\theta)$.

Since $f\in L^1(\R_+; E_{\theta})$, we get by the above conclusion
for ${T<\infty}$ and the representation formula \eq{2.1} for the
mild solution $u$ that the unique solution $u$ to \eq{2.8} satisfies
$u_t, Au\in L^1_{{\rm loc}}([0,\infty); E_{\theta}).$ Moreover,
since $\|u_t, Au\|_{L^1(0,T; E_{\theta})}$ can be estimated
independently of $T$, it follows that $u_t, Au\in L^1(\R_+;
E_{\theta})$ and the estimate \eq{2.20} holds true.
\hfill \qed\\[2ex]
{\bf Proof of Corollary \ref{C1.2}:} (i) For $\al\in (k,k+1)$,
$k\in\Z$, we have $E_\al=(E_k, E_{k+1})_{\al-k,1}$  and
$A_\al\in {\cal H}(E_{k+1},E_k)$, $\|e^{-tA_\al}\|_{{\cal
L}(E_k)}\leq M$ by \eq{2.6}.
Thus the assertion follows from Theorem \ref{T1.1} by considering the problem
in the underlying Banach space $E_k$.

 (ii) If $\al\in \Z$, then $A_{\al-1/2}\in
{\cal H}(E_{\al+1/2},E_{\al-1/2})$ and
$\|e^{-tA_{\al-1/2}}\|_{{\cal L}(E_{\al-1/2})}\leq M$ by \eq{2.6}.
Hence, $A^\bullet_\al$,  the $E^\bullet_\al$-realization of
$A_{\al-1/2}$ where $E^\bullet_\al:=(E_{\al-1/2},
E_{\al+1/2})_{1/2,1}$, has maximal $L^1$-regularity in
$E^\bullet_\al$ by Theorem \ref{T1.1}. \hfill\qed

\section{Maximal $L^1$-regularity of Dirichlet-Laplacian and Stokes operator in
Besov spaces $B^s_{q,1}$}
Important applications of the abstract theory derived in this paper
concern, in particular, maximal $L^1$-regularity of the
Dirichlet-Laplacian and the Stokes operator in inhomogeneous Besov
spaces.
Let $L^q(\Omega), W^{1,q}(\Omega)$ and $B^s_{q,r}(\Omega), 1\leq q,r\leq \infty,
s\in\R,$ denote the usual Lebesgue, Sobolev and Besov spaces,
respectively, on a domain $\Omega\subset \R^n$. If $\Omega\neq \R^n$,
the Besov space $B^s_{q,r}(\Omega)$ is defined by restriction of
tempered distributions in $B^s_{q,r}(\R^n)$ to $\Omega$;
its norm is defined by the quotient norm, see \cite{Am00}, \cite{Tri83}.

\subsection{Dirichlet-Laplacian}
Let $\Om\subset \R^n, n\in\N,$ be a domain with uniform
$C^2$-boundary $\pa\Om$. The Dirichlet-Laplacian $-\Da_\Om$ is
defined by
 $-\Da_\Om u:=-\Da u$ for $u$ in
$$ {\cal D}(-\Da_\Om):=\{u\in W^{2,q}(\Om): u|_{\pa\Om}=0\} = W^{2,q}(\Om) \cap W^{1,q}_0(\Om),\quad 1<q<\infty.$$
 It is well known that $-\Da_\Om$ generates an analytic
 semigroup in $L^q(\Om)$.
Hence, the $q$-dependent interpolation and extrapolation scale  $\{(E_{\al,1},
A_{\al,1}):\al\in\R\}$ generated by $(E,A):=(L^q(\Om), -\Da_\Om)$
and the real interpolation functor
$\{(\cdot,\cdot)_{\theta,1}:\theta\in (0,1)\}$ is well-defined. It
is known (\cite[Theorem 2.2, Proposition 2.4]{Am00}) that
$E_{\al,1} = B^{2\al}_{q,1,0}(\Om)$ for $\al\in (-1+1/2q,0)\cup
(0,1)$, where
\begin{equation}
\label{E4.1}
B^{s}_{q,1,0}(\Om)=\left\{
\begin{array}{cl}
\{u\in B^{s}_{q,1}(\Om): u|_{\pa\Om}=0\}, & s\in (1/q,2)\ek \{u\in
B^{s}_{q,1}(\R^n): {\rm supp}\, u\subset \bar\Om\}, & s=1/q\ek
 B^{s}_{q,1}(\Om),  & s\in (-2+1/q,1/q)\setminus\{0\},
  \end{array}
\right.
\end{equation}
and, of course, $B^{s}_{q,1,0}(\Om)=B^{s}_{q,1}(\R^n)$ if
$\Om=\R^n$. In general, $B^{s}_{q,1,0}(\Om)$ for $-2<s<0$ is defined
by
\begin{equation}
\label{E4.1n} B^{s}_{q,1,0}(\Om) = \big( (W^{2,q'}(\Om) \cap
W^{1,q'}_0(\Om))',L^q(\Om) \big)_{1+s/2,1},\; -2<s<0, \; q'=q/(q-1);
\end{equation}
however, if $\Om\neq \R^n$, the third characterization in \eq{4.1}
holds for $s<0$ only when $s\in (-2+1/q,0)$.

By the reiteration theorem yielding the identity
$(E_{-1/2,1},E_{1/2,1})_{1/2\pm\al,1}=E_{\pm\al,1}$ for sufficiently
small $\al>0$ ({\em cf.} \cite[Lemma 1.1]{Am00}), it also follows that
\begin{equation}
 \label{E4.2}
\begin{array}{rcl}
 E^\bullet_{0,1}&:= & (E_{-1/2,1},E_{1/2,1})_{1/2,1}\ek
  & = & \big((E_{-1/2,1},E_{1/2,1})_{1/2-\al,1},
(E_{-1/2,1},E_{1/2,1})_{1/2+\al,1}\big)_{1/2,1}\ek
 & = & (E_{-\al,1}, E_{\al,1})_{1/2,1}\quad (0<\al<1/2q)\ek
 & = & (B^{-2\al}_{q,1}(\Om),B^{2\al}_{q,1}(\Om))_{1/2,1}=B^0_{q,1}(\Om)=:B^0_{q,1,0}(\Om).
 \end{array}
\end{equation}

Thus, Theorem \ref{T1.1}, Corollary \ref{C1.2} result in the
following proposition where $-\Da_\Om$ denotes the
Dirichlet-Laplacian extended or restricted to $B^s_{q,1,0}(\Om)$,
$s\in (-2,2)$, see \eq{4.1}, \eq{4.1n} and \eq{4.2}.

\begin{propo}
\label{P4.1} {Let $\Om\in\R^n, n\in\N,$ be a domain with uniform
$C^2$-boundary such that semigroup $e^{t\Da_\Om}$ is bounded. Then
the Dirichlet-Laplacian $-\Da_\Om$ has maximal $L^1$-regularity in
$B^s_{q,1,0}(\Om)$, $q\in (1,\infty)$, $s\in (-2,2)$.}
\end{propo}
\begin{rem}{\em
If $\Om\neq \R^n$, the space $B^{s}_{q,1,0}(\Om)$ for $s\in
(-2, -2+1/q]$ in \eq{4.1n} can not be identified with a subspace of
${\mathcal D}'(\Om)$, the space of Schwartz distributions on $\Om$.
In fact, with the notation of \cite{Am03JNMP}, in this case
$B^{s}_{q,1,0}(\Om)=\big(\mathring B^{-s}_{q',\infty,0}(\Om)\big)'$
where
\begin{align}
\mathring B^{-s}_{q',\infty,0}(\Om) := &\, \big( W^{2,q'}(\Om) \cap
W^{1,q'}_0(\Om),L^{q'}(\Om) \big)_{1+s/2,\infty}^0\nonumber \ek
= &\, \{u\in \mathring B^{-s}_{q',\infty}(\Om): u|_{\pa\Om}=0\},\label{E4.4n}
\intertext{and}
\mathring B^{-s}_{q',\infty}(\Om):= &\, \big( W^{2,q'}(\Om), L^{q'}(\Om)
\big)_{1+s/2,\infty}^0,\; s\in (-2, -2+1/q],\nonumber
\end{align}
is the closure of $W^{-s,q'}(\Om)$ in the Besov space
$B_{q',\infty}^{-s}(\Om)$ ({\em cf.} \cite[(2.19), (2.20)]{Am00}
with $\mathring B^{-s}_{q',\infty} := n^{-s}_{q',\infty}$) yielding
$\big(\mathring B^{-s}_{q',\infty}(\Om)\big)'= B^s_{q,1}(\Om)$.
Moreover, for $s\in (-2, -2+1/q]$ the set ${\cal D}(\Om)$ is dense
in neither $\mathring B^{-s}_{q',\infty,0}(\Om)$ nor  $\mathring
B^{-s}_{q',\infty}(\Om)$, and neither $B^s_{q,1}(\Om)$ nor
$B^s_{q,1,0}(\Om)$ are subspaces of distribution spaces.
}
\end{rem}

\subsection{Stokes operator}

Maximal regularity of the Stokes operator is a crucial tool in the study of the
Navier-Stokes equations. Below, we briefly mention
how maximal $L^1$-regularity of the Stokes operator in solenoidal subspaces of  inhomogeneous
Besov spaces of $B^s_{q,1}$-type is implied by the abstract theory
of this paper.

 Let $\Om\in\R^n, n\geq 2,$ be a domain with uniform
 $C^2$-boundary $\pa\Om$. Let $1<q<\infty$ and $L^q_\si(\Om)$ and $W^{1,q}_{0,\si}(\Om)$
 be the closure of the set
$C^\infty_{0,\si}(\Om)=\{u\in C^\infty_{0}(\Om)^n:\div u=0\}$ in
$L^q(\Om)^n$ and $W^{1,q}(\Om)^n$, respectively. Assume that
 \begin{equation}
 \label{E4.5}
L^q_\si(\Om):=\{u\in L^q(\Om)^n:\div u=0,\;u\cdot {\bf
n}|_{\pa\Om}=0\},
 \end{equation}
  where ${\bf n}$ is the outward normal vector
at $\pa\Om$, and the {\it Helmholtz decomposition}
 \begin{equation}
 \label{E4.6}
 L^q(\Om)^n=L^q_\si(\Om)\oplus G_q(\Om),\quad
G_q(\Om)=\{v\in L^q(\Om)^n: \exists \pi\in L^1_{\rm loc}(\overline\Om): v=\na \pi \}
 \end{equation}
holds algebraically and topologically. Let ${\mathbb P} = {\mathbb P_q}$ be the {\it Helmholtz
projection} from $L^q(\Om)^n$ onto $L^q_\si(\Om)$. Then the Stokes
operator ${\mathbb A} = {\mathbb A_q}$ may be defined by
 \begin{equation}
 \label{E4.8}
 {\mathbb A}u:=-{\mathbb P}\Da u\; \text{ for }\;
u\in {\cal D}({\mathbb A}):=
  W^{2,q}(\Om)^n \cap W^{1,q}_{0,\si}(\Om),\; 1<q<\infty.
 \end{equation}
We assume that
\begin{equation}
 \label{E4.7}
{\mathbb A}\in {\mathcal H}({\mathbb E}_1, {\mathbb E}), \quad
\|e^{-t{\mathbb A}}\|_{{\cal L}({\mathbb E})}\leq Me^{\om
t}\quad(\exists M>0,\exists\om\geq 0,\forall t\geq 0),
 \end{equation}
where  ${\mathbb E}_1:={\cal D}({\mathbb A})$ is endowed with its
graph norm and ${\mathbb E}:=L^q_\si(\Om)$.

It is well known that the assumptions \eq{4.5}, \eq{4.6} and
\eq{4.7} are satisfied for many kinds of {\it standard} domains,
such as $\Omega=\R^n, \R^n_+$, bounded and exterior domains,
infinite layers and cylinders (with $\om=0$ in \eq{4.7}) and
aperture domains, and compact perturbations thereof (with $\om\neq
0$, in general), see Introduction of \cite{FR07} and the references
therein. However, there are smooth unbounded domains for which the
Helmholtz projection does not exist, see \cite{FKS05}.

Starting from $(\mathbb{E}, {\mathbb A})$ and the real interpolation functor
$\{(\cdot,\cdot)_{\theta,1}:\theta\in (0,1)\}$, the interpolation
and extrapolation scale  $\{({\mathbb E}_{\al,1}, {\mathbb
A}_{\al,1}):\al\in\R\}$ is generated. There exists an explicit
representation of ${\mathbb E}_{\al,1}$ for $|\al|\leq 1$, see
\cite[Theorem 3.4, Remark 3.7]{Am00}.
 More precisely, using the notation from \cite{Am03JNMP},
 we have
${\mathbb E}_{\al,1} = \mathbb B^{2\al}_{q,1}(\Om)$ with
\begin{equation}
\label{E4.3}  \mathbb B^{s}_{q,1}(\Om):=\left\{
\begin{array}{cl}
\{u\in B^{s}_{q,1}(\Om)^n: \div u=0,\; u|_{\pa\Om}=0\}, & s\in
(1/q,2)\ek
 \{u\in B^{s}_{q,1}(\R^n)^n: \div u=0,\;  {\rm supp}\, u\subset \bar\Om\}, & s=1/q\ek
 \{u\in B^{s}_{q,1}(\Om)^n: \div u=0,\;u\cdot {\bf n}|_{\pa\Om}=0\},  & s\in
 (0,1/q) \ek
\big(\mathcal D(\mathbb A_{q'})',L^q_\sigma(\Om)\big)_{1+s/2,1}, & s\in(-2,0).
 \end{array}
\right.
\end{equation}

Moreover, for ${\mathbb E}^\bullet_{0,1}:=({\mathbb
E}_{-1/2,1},{\mathbb E}_{1/2,1})_{1/2,1}$ we have
 \begin{equation}
 \label{E4.4}
 \mathbb E^\bullet_{0,1}=\{u\in B^0_{q,1}(\Om)^n: \div u=0,\; u\cdot {\bf
 n}|_{\pa\Om}=0\}=: \mathbb B^{0}_{q,1}(\Om),
 \end{equation}
 which can be proved exactly in the same way as \cite[Lemma 4.3]{RiZhZh15} just replacing
 $E_{s/2,\infty}(q)$, $B^s_{q,\infty}$ by ${\mathbb E}_{s/2,1}$,
$B^s_{q,1}$, respectively.

Finally, applying Theorem \ref{T1.1} and Corollary \ref{C1.2}, we
obtain the following result.

\begin{propo}
\label{P4.2} Let $\Om\in\R^n, n\geq 2,$ be ${\mathbb R}^n, {\mathbb
R}^n_+$, infinite layers, or bounded, exterior domains, or cylinders
with uniform $C^2$-boundaries. Moreover, let the Stokes operator
${\mathbb A}$ be given as in \eq{4.8}. Then ${\mathbb A}$ has
maximal $L^1$-regularity in $\mathbb B^s_{q,1}(\Om)$, $q\in
(1,\infty)$, $|s|<2$.
\end{propo}
\begin{rem}{\rm
\label{R4.4} (1)
If $\Om=\R^n$, it follows by \cite[Remark 3.7
(a)]{Am00} that for $1<q<\infty$ and $0<|s|<2$ there holds
\begin{equation}\label{E4.10a}
\mathbb B^s_{q,1}(\R^n) =\{u\in B^{s}_{q,1}(\R^n)^n:\div u=0\}.
\end{equation}

(2) If $\Om\neq \R^n$ and $s<0$, an explicit characterization
of $\mathbb B^s_{q,1}(\Om)$ is, in general, difficult to find, see
\cite{Am00}. If $s\in (-2,-2+1/q]$, then $(E_{s/2,1})^n$ is not a
subspace of distributions on $\Om$.

For $s\in (-2+1/q,0)$ a characterization of $\mathbb B^s_{q,1}(\Omega)$ similar to \eqref{E4.10a} is possible provided there exists a Helmholtz decomposition of
$(E_{s/2,1})^n=B^s_{q,1}(\Om)^n$ in the form
 \begin{equation}
 \label{E4.10}
B^s_{q,1}(\Om)^n=M\oplus N
\end{equation}
 with
 $M:=\{u\in B^{s}_{q,1}(\Om)^n:\div u=0\}$ and $N:=\{\na p\in (B^{s}_{q,1}(\Om)^n:
 \;\exists p\in {\cal D}'(\Om)\}$.
 From the
duality property \eq{3.2} one infers that $\mathbb B^s_{q,1}(\Om)$
is the dual of $\mathring{\mathbb B}^{-s}_{q',\infty}(\Om)$, where
$\mathring{\mathbb B}^{-s}_{q',\infty}(\Om)$ is the completion of
${\mathcal D}(\mathbb A_{q'})$ in $B_{q',\infty}^{-s}(\Om)^n$. Note
that $\mathring{\mathbb B}^{-s}_{q',\infty}(\Om)\subset \mathring{
B}^{-s}_{q',\infty}(\Om)^n$ and $({\mathring
B}^{-s}_{q',\infty}(\Om))'=B^s_{q,1}(\Om)$, see \S 4.1. Then, by
Hahn-Banach's theorem,
 \begin{equation}
 \label{E4.12}
\mathbb B^s_{q,1}(\Om)\cong B^{s}_{q,1}(\Om)^n/\tilde{N},\;\;s\in
(-2+1/q,0),
 \end{equation}
  where $"\cong"$ means isometric isomorphism, and
$$ \tilde{N}=\big\{u\in B^s_{q,1}(\Om)^n:
 \lan u,\vp\ran_{B^s_{q,1}, \mathring{ B}^{-s}_{q',\infty}}=0\quad \forall
 \vp\in \mathring{\mathbb B}^{-s}_{q',\infty}(\Om)\big\}. $$
 By de Rham's lemma we have $\tilde{N}=N$. Then, thanks to \eq{4.10} and \eq{4.12},
$$ \mathbb B^s_{q,1}(\Om) \cong M=\{u\in B^{s}_{q,1}(\Om)^n:\div u=0\}. $$
}
\end{rem}

\renewcommand{\baselinestretch}{1.0} 

\end{document}